\documentclass[11pt]{article}
\usepackage{amssymb, amsmath, amsthm, graphicx, hyperref}
\usepackage[left=1in,top=1in,right=1in]{geometry}
\date{}
\allowdisplaybreaks

\theoremstyle{definition}
\newtheorem{theorem}{Theorem}[section]

\newtheorem{problem}[theorem]{Problem}

\newtheorem{corollary}[theorem]{Corollary}

\counterwithin{equation}{section}

\title{3-uniform monotone paths and multicolor Ramsey numbers}
\author{Andrew Suk\thanks{Department of Mathematics, University of California San Diego, La Jolla, CA, 92093 USA. Research supported by NSF Awards DMS-1952786 and DMS-2246847. Email: {\tt asuk@ucsd.edu}.} \and Ji Zeng\thanks{Alfr\'ed R\'enyi Institute of Mathematics, Budapest, Hungary. Supported by European Research Council grants No. 882971 and No. 101054936. Email: {\tt jzeng@ucsd.edu}.}}

\begin{document}

\maketitle

\begin{abstract}
    The \textit{monotone path} $P_{n+2}$ is an ordered 3-uniform hypergraph whose vertex set has size $n+2$ and edge set consists of all consecutive triples. In this note, we consider the collection $\mathcal{J}_n$ of ordered 3-uniform hypergraphs named \textit{monotone paths with $n$ jumps}, and we prove the following relation
    \begin{equation*}
        r(3;n) \leq R(P_{n+2},\mathcal{J}_n) \leq 4^n \cdot r(3;n),
    \end{equation*} where $r(3;n)$ is the multicolor Ramsey number for triangles and $R(P_{n+2},\mathcal{J}_n)$ is the hypergraph Ramsey number for $P_{n+2}$ versus any member of $\mathcal{J}_n$. In particular, whether $r(3;n)$ is exponential, which is a very old problem of Erd\H{o}s, is equivalent to whether $R(P_{n+2},\mathcal{J}_n)$ is exponential.
\end{abstract}

\section{Introduction}

The \textit{multicolor Ramsey number} $r(m;n)$ is the smallest number $N$ such that any $n$-coloring of the edges of the complete graph on $N$ vertices contains a monochromatic complete subgraph on $m$ vertices. The existence of $r(m;n)$ follows from the celebrated theorem of Ramsey \cite{ramsey} from 1930, and in the special case when $m = 3$, Issai Schur proved the existence of $r(3;n)$ in 1916 in his work related to Fermat’s Last Theorem \cite{schur}.  In particular, he showed that

$$2^{\Omega(n)} < r(3;n) < O(n!).$$

\noindent  While the upper bound has remained unchanged over the last 100 years, the best known lower bound $r(3;n) > 3.199^n$ can be obtained by combining the results of Fredricksen and Sweet \cite{FS} and Abbot and Moser \cite{AB} (see also \cite{Xi,conlon2021lower}). It is a major open problem in Ramsey theory to close the gap between the lower and upper bounds for $r(3;n)$. Erd\H{o}s even offered cash prizes for solutions to the following problems (see \cite{chung1998erdos}).

\begin{problem}[\$250]
    Determine $\lim_{n \to \infty} \left( r(3;n) \right)^{1/n}$.
\end{problem}

It was shown by Chung~\cite{chung1973ramsey} that $r(3; n)$ is supermultiplicative, so the above limit exists.

\begin{problem}[\$100]\label{prob}
    Is $\lim_{n \to \infty} \left( r(3;n) \right)^{1/n}$ finite or not?
\end{problem} The goal of this note is to show that Problem~\ref{prob} is equivalent to another Ramsey-type problem for ordered 3-uniform hypergraphs, which we describe below.

A \textit{monotone path} $P_n$ is a vertex-ordered 3-uniform hypergraph whose vertex set is $\{1,2,\ldots, n\}$ and edge set consists of all consecutive triples $(i,i+1,i+2)$. The Ramsey number $R(P_n,P_n)$ is the smallest number $N$ such that any red-blue coloring of the complete (vertex-)ordered $3$-uniform hypergraph on $N$ vertices, denoted as $K^{(3)}_N$, contains a monochromatic copy of $P_n$. A famous result by Erd\H{o}s and Szekeres~\cite{erdos1935combinatorial} states that \begin{equation*}
    R(P_n,P_n) = \binom{2n-4}{n-2}+1 = 4^{n - o(n)}.
\end{equation*}

We consider the collection $\mathcal{J}_n$ of ordered 3-uniform hypergraphs $H$ such that $H$ contains all consecutive triples as its edges and a vertex subset $J\subset V(H)$ of size $n$ satisfying the following conditions:
\begin{enumerate}
    \item[(0)]  The subset $J$ does not contain consecutive pairs, the first vertex, or the last vertex of $H$.
    \item[(1)]  For each $v \in J$, the triples $(v-2,v-1,v+1)$ and $(v-1,v+1,v+2)$ are edges in $H$.
    \item[(2)]  If both $v-1$ and $v+1$ are in $J$, then  $(v-2,v,v+2)$ is also an edge in $H$.
\end{enumerate} Here, $v+1$ denotes the vertex that immediately succeeds $v$ in $H$, and the meanings are similar for $v+2$, $v-1$, etc.

We refer to members of $\mathcal{J}_n$ as \textit{monotone paths with $n$ jumps} and their corresponding $J$ as \textit{jumps}. Intuitively, $H \in \mathcal{J}_n$ still contains a monotone path from start to end if any jump is ``skipped''. Let $R(P_{n+2},\mathcal{J}_n)$ be the smallest number $N$ such that any red-blue coloring of $K^{(3)}_N$ contains a red $P_{n+2}$ or a blue member of $\mathcal{J}_n$. Our main result is the following relation.
\begin{theorem}\label{main}
    We have $r(3;n) \leq R(P_{n+2},\mathcal{J}_n) \leq 4^n \cdot r(3;n)$ for all positive integer $n$.
\end{theorem}

In particular, whether $r(3;n)$ is exponential in $\Theta(n)$ equates with whether $R(P_n,\mathcal{J}_n)$ is exponential in $\Theta(n)$. We remark that monotone paths with $n$ jumps are quite similar to monotone paths, say, in terms of edge density. For example, let $P^t_n$ denote the $t$-th power path of $P_n$, precisely, the ordered 3-uniform hypergraph on vertices $\{1,2,\ldots, n\}$ such that every $t$ consecutive vertices form a $K_t^{(3)}$. In particular, $P_n^3  = P_n$. Notice that $P_{3n}^4$ contains a monotone path with $n$ jumps, where the jumps occur at $3i-1$.
\begin{corollary}
    We have  $r(3;n) \leq R(P_{n+2},P^4_{3n})$ for all positive integer $n$.
\end{corollary}

Recently, Ramsey-type results for power paths have been studied in the ordered graph setting, and we refer the interested reader to \cite{Su,MS}. In view of the Erd\H{o}s--Szekeres theorem stated above, as well as the similarity between $P_n$ and $P^4_{n}$, we suspect that the Ramsey number $r(3;n)$ is exponential in $\Theta(n)$.

\section{Proof}

\begin{proof}[Proof of lower bound]

    Let $N = r(3;n) - 1$ and $K^{(3)}_N$ be the complete ordered 3-uniform hypergraph with vertex set $[N] = \{1,2,\dots,N\}$. By the definition of $r(3;n)$, there exists an $n$-coloring $\chi$ of the pairs of $[N]$ that avoids monochromatic triangles. We color each edge $(u,v,w)$ of $K^{(3)}_N$ red if $\chi(u,v) < \chi(v,w)$, otherwise we color it blue. It suffices to show this colored $K^{(3)}_N$ does not contain a red $P_{n+2}$ or a blue member of $\mathcal{J}_n$.

For sake of contradiction, suppose inside $K^{(3)}_N$ there is a red $P_{n+2}$ with vertices $v_1 < \dots< v_{n+2}$, then we have\begin{equation*}
        1\leq \chi(v_1,v_2) < \chi(v_2,v_3) < \dots < \chi(v_{n+1}, v_{n+2}) \leq n.
    \end{equation*} This is a contradiction since there are not as many distinct integers in $[n]$.

    For each $H \in \mathcal{J}_n$ with jumps $J\subset V(H)$, we associate an ordered (2-uniform) graph $G_H$ whose vertex set is $V(H)$, and edge set consists of all consecutive pairs in $V(H)$, and all pairs $(v-1,v+1)$ with $v\in J$. In what follows, we prove by induction on $n$ that: if $\chi$ is an $n$-coloring of the edges of $G_H$ associated to $H \in \mathcal{J}_n$, where $\chi(u,v) \geq \chi(v,w)$ for all $(u,v,w) \in E(H)$, then there is a monochromatic triangle in $G_H$ under $\chi$. Given this claim, $K^{(3)}_N$ cannot contain a blue member of $\mathcal{J}_n$, hence concluding the proof.

    When $n=1$, we pick $v \in J \subset V(H)$. Condition (0) of $\mathcal{J}_n$ guarantees the existence of $v-1$ and $v+1$ in $H$. Clearly, $\{v-1,v,v+1\}$ form a monochromatic triangle in $G_H$ under $\chi$. This proves the base case.
    
    For the inductive step, when $n>1$, let $c$ be the smallest color in $G_H$ under $\chi$, and $w$ be the largest vertex of $J$. If the color $c$ only appears among the edges $(w-1,w)$, $(w,w+1)$, and $(w-1,w+1)$, we consider the sub-hypergraph $H'\subset H$ induced on vertices that appear before $w$. It is easy to check that $H' \in \mathcal{J}_{n-1}$ and $G_{H'}$ only receives $n-1$ colors under $\chi$ (because $c$ is not used). Then we can apply the inductive hypothesis to conclude this case.

    If $c$ appears in the edge $(w-2 , w-1)$, we argue $\{w-1,w,w+1\}$ forms a monochromatic triangle with color $c$. Indeed, by condition (1) of $\mathcal{J}_n$, $(w-2,w-1,w+1)\in E(H)$, so $\chi(w-2,w-1) \geq \chi(w-1,w+1)$. Then the minimality of $c$ implies $\chi(w-1,w+1) = c$. Similarly, since $E(H)$ includes all consecutive triples in $H$, we can argue $\chi(w-1,w) = \chi(w,w+1) = c$.

    If $c$ appears in the edge $(v, w-1)$ with $v \neq w-2$, by the definition of $G_H$, we must have $v=w-3$ and $w-2\in J$. In this case, Condition (2) of $\mathcal{J}_n$ implies that $(w-3, w-1 ,w+1)\in E(H)$. Other conditions of $\mathcal{J}_n$ can be used to show $(w-3,w-1,w),(w-1,w,w+1) \in E(H)$. Then we similarly argue that $\{w-1,w,w+1\}$ forms a monochromatic triangle with color $c$.

    If $c$ appears in an edge $(u,v)$ with $v < w - 1$, we argue that $\chi(v,v+1) = c$ as well. Indeed, $(u,v,v+1) \in E(H)$ follows from conditions of $\mathcal{J}_n$, then we can apply the minimality of $c$. Iterating this process, we can reach $\chi(w-2,w-1) = c$, and this forces a monochromatic triangle $\{w-1,w,w+1\}$ as argued above, finishing the inductive argument.
\end{proof}

\begin{proof}[Proof of upper bound]

    We consider a specific member $I_n \in \mathcal{J}_n$ whose vertex set is $[2n+1]$, and edge set consists of all consecutive triples $(i,i+1,i+2)$ as well as all triples in the form of $(2i-1,2i+1,2i+2)$, $(2i, 2i+1, 2i+3)$, and $(2i-1, 2i+1, 2i+3)$. See Figure~\ref{fig:I2} for an example. It is easy to check that $\{2i: i \in [n]\}$ are jumps of $I_n$. Assume $K^{(3)}_N$ is endowed with a red-blue coloring on its edges (triples), and it does not contain a red $P_{n}$ or a blue $I_n$. It suffices to show $N \leq 4^n \cdot r(3;n)$.

    \begin{figure}[ht]
        \centering
        \includegraphics[scale=1]{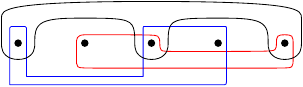}
        \caption{Edges of $I_2$ that are not consecutive triples.}
        \label{fig:I2}
    \end{figure}

    For every two vertices $u<v$ of $K^{(3)}_N$, we define \begin{equation*}
        \alpha(u,v) = 1 + \text{number of edges in the longest red monotone path ending at $(u,v)$}.
    \end{equation*} We define $\beta(u,v)$ to be the largest $\ell$ such that there exists $v_1 < v_2 < \dots < v_{2\ell - 1}$ in $[N]$ such that $v_{2\ell-2} = u$, $v_{2\ell-1} = v$ and for all possible $i$,\begin{equation*}
        \alpha(v_{2i-1}, v_{2i}) =  \alpha(v_{2i}, v_{2i+1}) =  \alpha(v_{2i-1}, v_{2i+1}) \geq \alpha(v_{2i+1}, v_{2i+2}) =  \alpha(v_{2i+2}, v_{2i+3}) =  \alpha(v_{2i+1}, v_{2i+3}).
    \end{equation*} If such $\ell$ does not exist, we define $\beta(u,v) = 1$. Notice that $\alpha(u,v) \geq \alpha(v,w)$ implies the triple $(u,v,w)$ is blue. Using this fact we can check that $\beta(u,v) > n$ would force a blue $I_n$ inside $K^{(3)}_N$. Therefore, we have $1\leq \alpha(u,v),\beta(u,v) \leq n$ for all $(u,v)$.

    For each vertex $v$, we associate a set $D(v) \subset [n]^2$ defined as follows
    \begin{equation*}
        D(v) = \{(a,b):~\exists u<v,~ a\leq \alpha(u,v),~b\leq \beta(u,v)\}.
    \end{equation*} Notice that $D(v)$ is ``downward-closed'', that is, $(a,b) \in D(v)$ implies $(a',b') \in D(v)$ for $a' \leq a$ and $b' \leq b$. Using the bijection between such subsets and integer partitions (see \cite{moshkovitz2014ramsey} and \cite{chvatal1971some}), we can upper bound the total number of ``downward-closed'' subsets by $\binom{2n}{n} \leq 4^n$. Suppose $N > 4^n \cdot r(3;n)$ for a contradiction, there must be $r(3;n)$ vertices, denoted as $S$, associated with the same ``downward-closed'' set $\Delta$.
    
    Now regard $\alpha$ as an $n$-coloring of pairs of $S$. By definition of $r(3;n)$, there is a monochromatic triangle $\{u,v,w\}$ under $\alpha$. Since $(\alpha(v,w),\beta(v,w)) \in D(w)$ and $D(w) = \Delta =D(u)$, there exists $t < u$ such that $\alpha(v,w) \leq \alpha(t,u)$ and $\beta(v,w) \leq \beta(t,u)$. By definition of $\beta(t,u)$, there exists a sequence $v_1 < v_2 < \dots < v_{2\ell - 1}$ with $\ell = \beta(t,u)$ such that $v_{2\ell-2} = t$, $v_{2\ell-1} = u$ and \begin{equation*}
        \alpha(v_{2i-1}, v_{2i}) =  \alpha(v_{2i}, v_{2i+1}) =  \alpha(v_{2i-1}, v_{2i+1}) \geq \alpha(v_{2i+1}, v_{2i+2}) =  \alpha(v_{2i+2}, v_{2i+3}) =  \alpha(v_{2i+1}, v_{2i+3}).
    \end{equation*} (We remark that the monochromatic triangle $\{u,v,w\}$ already guarantees $\beta(v,w) \geq 2$ hence the case $\beta(t,u) = 1$ is non-existent.) Notice that $\alpha(v_{2\ell-2},v_{2\ell-1}) \geq \alpha(v,w) = \alpha(u,v) = \alpha(u,w)$. As a result, the sequence $v_1 < \dots < v_{2\ell - 1} < v < w$ will imply $\beta(v,w) \geq \ell + 1 > \beta(v,w)$, which is a contradiction. This concludes our proof.
\end{proof}

\section{Concluding Remarks}

Our proof of lower bound is inspired by a result of Mubayi and Suk~\cite{mubayi2017off}. Our proof of upper bound is inspired by the works of Moshkovitz--Shapira~\cite{moshkovitz2014ramsey} and Chv{\'a}tal--Koml{\'o}s~\cite{chvatal1971some}.

It is obvious from the proof (but not the statement) that if one desires to prove $r(3;n)$ is super-exponential in $n$, it suffices to show that $R(P_n,I_n)$ is super-exponential in $n$ for the specific $I_n$ as defined above.  On the other hand, to prove $r(3;n)$ is exponential in $\Theta(n)$, it may be unnecessary to show $R(P_{n},\mathcal{J}_n)$ is exponential in $\Theta(n)$. One possible approach is to find families of non-increasing triples inside $n$-colored (unordered and 2-uniform) complete graphs that form monotone paths with $n$ jumps after some suitable orderings.

It is plausible that our method could be used to prove $r(m;n) \leq R(P_{n}, H_n)$ where $H_n$ is an ordered 3-uniform hypergraph whose vertex set has size linear in $n$ and edge-density arbitrarily close to that of the monotone path. However, to exactly describe all such $H_n$ could be a tedious task. Nevertheless, the proof of Theorem \ref{main} can be easily modified to show the following.

\begin{theorem}
        We have $r(m;n) \leq R(P_{n+2}, P_{mn}^{m + 1})$ for all integers $m,n\geq 3$.
\end{theorem}

\medskip

\noindent {\bf Acknowledgement.} We wish to thank G\'abor Tardos for several helpful comments.

\bibliographystyle{abbrv}
{\footnotesize \bibliography{main}}

\end{document}